\documentclass[final]{article}

\usepackage{graphicx}
\usepackage{color}
\usepackage{amsmath,amsthm,amssymb,amsfonts}
\usepackage{mathrsfs}
\usepackage{url}
\usepackage{geometry} 
\usepackage{enumerate}

\newcounter{cst}
\def \ctel#1{C_{\refstepcounter{cst}\label{#1}\thecst}} 
\def \cter#1{C_{\ref{#1}}} 

\newtheorem{theorem}{Theorem}[section] 
\newtheorem{remark}[theorem]{Remark}
\newtheorem{lemma}[theorem]{Lemma}

\newenvironment{sketch}{
  \proof}{\endproof}

\def\bhyp#1{\begin{equation}\label{#1}\begin{array}{c}}
\def\ehyp{\end{array}\end{equation}}

\newcommand{\RR}{{\mathbb R}}
\newcommand{\NN}{{\mathbb N}}

\renewcommand{\O}{\Omega}

\newcommand{\cF}{\mathcal{F}}
\newcommand{\cD}{\mathcal{D}}
\newcommand{\cS}{\mathcal{S}}
\newcommand{\cI}{\mathcal{I}}

\newcommand{\ch}{\mathbf{1}}
\newcommand{\U}{\mathbf{u}}

\newcommand{\Ubar}{\overline{\U}}
\newcommand{\cbar}{\overline{c}}
\newcommand{\pbar}{\overline{p}}

\newcommand{\weakto}{\rightharpoonup}

\newcommand{\weak}{\mbox{\rm-w}}

\newcommand{\edges}{\mathcal{E}}
\newcommand{\Kedges}{\edges_{K}}
\newcommand{\intedges}{\edges_{\rm{int}}}

\newcommand{\cells}{\mathcal{M}}
\newcommand{\normvec}{\mathbf{n}}
\newcommand{\mesh}{\cD}
\newcommand{\fluxes}{\cF_\mesh}
\newcommand{\unknowns}{X_\mesh}

\newcommand{\deltat}{\delta t^{(n-\frac{1}{2})}}

\newcommand{\dive}{\operatorname{div}} 
\newcommand{\ud}{\, \mathrm{d}} 
\newcommand{\norm}[1]{\left\lVert#1\right\rVert} 
\newcommand{\size}[1]{\operatorname{size}(#1)} 

\newcommand{\Leb}[3][ ]{L^{#2}(0, T; L^{#3}(\O)^{#1})} 

\newcommand{\D}[3][ ]{\mathbf{D}_{#1}(#2, #3)} 

\newcommand{\Dsq}[3][ ]{\mathbf{D}_{#1}^{1/2}(#2, #3)} 

\newcommand{\DD}{\mathbf{D}}

\title{Uniform temporal convergence of numerical schemes for incompressible
miscible displacement}
\author{Kyle S. Talbot\footnote{School of Mathematical Sciences, Monash University, 
Victoria 3800, Australia. \texttt{kyle.talbot@monash.edu}}}
\date{\today}

\begin{document}
\maketitle

\begin{abstract}
The Hybrid Mimetic Mixed (HMM) family of discretisations includes the Hybrid
Finite Volume method, the Mimetic Finite Difference method and the Mixed Finite
Volume method. This note demonstrates that HMM discretisations
of the equations describing the single-phase, miscible displacement through a 
porous medium of one incompressible fluid by another converge uniformly in time 
for the concentration variable.
\end{abstract}

\maketitle

\section{Introduction} \label{sec:intro}

The nonlinearly-coupled elliptic-parabolic system
\begin{equation}\left.
\begin{aligned}
&\dive(\Ubar) = q^I - q^P \quad\mbox{in $\O\times(0,T)$,}&\qquad
&\Ubar = -A(\cdot,\cbar)\nabla\pbar \quad\mbox{in $\O\times(0,T)$,}\\
&\int_\O\pbar(x,\cdot)\ud x = 0 \quad\mbox{on $(0,T)$,}&\qquad
&\Ubar\cdot\normvec=0 \quad\mbox{on $\partial\O\times(0,T)$,}
\end{aligned}\right\} \label{eq:elliptic}
\end{equation}
\begin{equation}\left.
\begin{aligned}
&\Phi\partial_{t}\cbar - \dive(\D{\cdot}{\Ubar}\nabla\cbar - \cbar\,\Ubar)
= \hat{c}q^I - \cbar q^{P}\quad\mbox{in $\O\times(0,T)$,}\\
&\cbar(\cdot,0) = c_0 \quad\mbox{in $\O$,}\\
&\D{\cdot}{\Ubar}\nabla\cbar\cdot\normvec = 0 \quad\mbox{on $\partial\O\times(0,T)$,}
\end{aligned}\right\} \label{eq:parabolic}
\end{equation}
describes the single-phase, miscible displacement through a porous medium of one incompressible
fluid by another, in the absence of gravity \cite{pe77}.
The unknowns are the pressure $\pbar$ of the fluid mixture, the Darcy velocity $\Ubar$
of the fluid mixture and the concentration $\cbar$ of the injected fluid in the medium.
We assume that the porous medium $\O$ is an open, bounded, convex polygonal subset 
of $\RR^{d}$, $d\geq2$, and that the displacement occurs over the time interval
$(0,T)$, $T>0$. The porosity $\Phi\in L^{\infty}(\O)$, and there is $\phi_{\ast}>0$
such that for a.e. $x\in\O$, $\phi_{\ast}\leq\Phi(x)\leq\phi_{\ast}^{-1}$. The
coefficient $A:\O\times\RR\to M_{d}(\RR)$ is a uniformly elliptic, bounded
Carath\'eodory function that combines the absolute permeability of the medium
and viscosity of the fluid mixture. The injected concentration $\hat{c}\in L^{\infty}(\O\times(0,T))$
satisfies $0\leq\hat{c}(x, t)\leq 1$ for a.e. $(x, t)\in\O\times(0, T)$,
and the initial concentration $c_0\in L^{\infty}(\O)$ satisfies $0\leq c_{0}(x)\leq 1$
for a.e. $x\in\O$. The injection well source terms $q^I\in\Leb{\infty}{2}$ and production well sink terms
$q^P\in\Leb{\infty}{r}$ (for some $r>2$) are nonnegative and satisfy 
$\int_\O q^I(x,t)\ud x=\int_\O q^P(x,t)\ud x$ for all $t\in(0,T)$. 
Finally, the diffusion-dispersion tensor $\DD:\O\times\RR^{d} \to M_{d}(\RR)$ is a Carath\'eodory
function with positive constants $\alpha_{\DD}$, $\Lambda_{\DD}$ such that for a.e. 
$x\in\O$ and all $\zeta,\xi\in\RR^{d}$,
\begin{equation} \label{hyp:D}
\DD(x,\zeta)\xi\cdot\xi \geq \alpha_{\DD}(1+|\zeta|)|\xi|^{2} \mbox{ and } 
|\DD(x,\zeta)|\leq \Lambda_{\DD}(1+|\zeta|).
\end{equation}
Our notion of weak solution to \eqref{eq:elliptic}--\eqref{eq:parabolic} is \cite[Definition 1.1]{cd07}.

\section{Preliminaries} \label{sec:scheme}

The HMM family of discretisations includes Hybrid Finite Volumes (HFV) \cite{egh10}, 
Mimetic Finite Differences \cite{bls05} and Mixed Finite Volumes (MFV) \cite{dey06}. 
These three families are equivalent \cite{degh10}.
An implementation of the HMM method amounts to a choice of any of these three discretisations.

The variable of most interest is the concentration $\cbar$. To obtain meaningful
approximations of this quantity, one must account for the (dominant) convective term
$\dive(\cbar\,\Ubar)$ in \eqref{eq:parabolic}. One therefore requires a good approximation 
of the normal flux of $\Ubar$ through the boundaries of control volumes of the mesh. 
We therefore follow \cite{cd07} by choosing the MFV framework to discretise \eqref{eq:elliptic}.
The HFV framework enables direct approximations of $\cbar$ and $\nabla\cbar$ without
first approximating the flux of $\D{\cdot}{\Ubar}\nabla\cbar$, which for our purposes
herein would be redundant. For brevity, we omit the details of the MFV scheme for
\eqref{eq:elliptic}; interested readers should consult \cite{cd07}.

With the exception of the convective term, we present the HFV discretisation of \eqref{eq:parabolic} 
in the notation of \emph{gradient schemes}, a framework introduced by \cite{degh13}. 
Incorporating the discretisation of convective terms into the gradient schemes 
framework is the subject of future work.

We adopt the same notion of admissible mesh $\mesh=(\cells,\edges)$ of $\O$ and corresponding notation as 
\cite{cd07}, with the following exceptions: we write $|K|$ 
for the $d$-dimensional measure of a control volume $K$, and $|\sigma|$ for the 
$(d-1)$-dimensional measure of an edge $\sigma$. 
For all $K\in\cells$ and $\sigma\in\Kedges$, we denote by $d_{K,\sigma}$ the 
Euclidean distance between $x_K$ and the hyperplane containing $\sigma$. 
Furthermore, we discretise the time interval by choosing a sequence $(t^{(n)})^{n=0,\ldots,N}$ such
that $0=t^{(0)}<t^{(1)}<\cdots<t^{(N)} = T$. For $n=1,\ldots, N$, set 
$\delta t^{(n-\frac{1}{2})}=t^{(n)} - t^{(n-1)}$ and
$\delta_\mesh =\max_{n=1,\ldots,N}\delta t^{(n-\frac{1}{2})}$. 

The space of discrete unknowns is 
$\unknowns := \left\{ v=((v_K)_{K\in\cells}, (v_\sigma)_{\sigma\in\edges}) 
: v_K\in\RR, v_\sigma\in\RR \right\}$. For $v\in\unknowns$ and $K\in\cells$, define 
the operator $\Pi_\mesh:\unknowns\to L^2(\O)$ by $\Pi_\mesh(v)=v_K$ on $K$.
We employ the same discrete gradient operator $\nabla_\mesh : \unknowns\to L^2(\O)^d$
as in \cite[Eq. (5.8)]{degh13}; see this reference for further details of the construction.
The norm on $v\in\unknowns$ is then
$\norm{v}_{\unknowns}:= \norm{\Pi_\mesh v}_{L^2(\O)} + \norm{\nabla_\mesh v}_{L^2(\O)^d}$.
To specify the initial condition in the scheme, we use a linear interpolation operator
$\cI_\mesh : L^2(\O)\to\unknowns$. The space of discrete fluxes is 
$\fluxes := \left\{ G=(G_{K,\sigma})_{K\in\cells,\sigma\in\Kedges} : G_{K,\sigma}\in\RR\right\}$.

In the scheme below, we consider sequences 
$(c^{(n)})_{n=0,\ldots,N}\subset\unknowns$ and $(F^{(n)})_{n=1,\ldots,N}\subset\fluxes$. 
For $n=1,\ldots,N$, $c^{(n)}_K$ approximates $\cbar$ on $K\times[t^{(n-1)}, t^{(n)})$, and 
$F^{(n)}_{K,\sigma}$ approximates the flux $-\int_\sigma\Ubar\cdot\normvec_{K,\sigma}\ud\gamma$
on $[t^{(n-1)}, t^{(n)})$. As in \cite{cd07}, we denote by $\U$ the piecewise-constant, 
cell-centred approximation of $\Ubar$ coming from the MFV scheme for \eqref{eq:elliptic}.

Following convention in the gradient schemes literature to date, we use the notation
$\Pi_\mesh$ and $\nabla_\mesh$ for functions dependent on both space and time. Thus
if $(v^{(n)})_{n=0,\ldots,N}\subset\unknowns$, for a.e. $x\in\O$ we set 
$\Pi_\mesh v(x,0)=\Pi_\mesh v^{(0)}(x)$, and for all $n=1,\ldots,N$, all 
$t\in[t^{(n-1)}, t^{(n)})$ and a.e. $x\in\O$, $\Pi_\mesh v(x,t)=\Pi_\mesh v^{(n)}(x)$, 
$\nabla_\mesh v(x,t)=\nabla_\mesh v^{(n)}(x)$ and
\begin{equation*}
\delta_\mesh v(t)=\delta_\mesh^{(n-\frac{1}{2})}v
:= \frac{v^{(n)}-v^{(n-1)}}{\delta t^{(n-\frac{1}{2})}}\in\unknowns.
\end{equation*}

The scheme for \eqref{eq:parabolic} is then: find sequences
$(c^{(n)})_{n=0,\ldots,N}\subset\unknowns$ and $(F^{(n)})_{n=1,\ldots,N}\subset\fluxes$ such that
\begin{equation}\left.
	\begin{gathered}
	c^{(0)}=\cI_\mesh c_0 \mbox{ and for all $\varphi=(\varphi^{(n)})_{n=1,\ldots,N}\subset\unknowns$,}\\
	\int_0^T\int_\O\left[\Phi(x)\Pi_\mesh\delta_\mesh c(x,t)\Pi_\mesh\varphi(x,t)
	+ \D{x}{\U(x,t)}\nabla_\mesh c(x,t)\cdot\nabla_\mesh\varphi(x,t)\right]\ud x\ud t\\
	+\sum_{n=1}^N\deltat \sum_{K\in\cells}\sum_{\substack{\sigma\in\Kedges\cap\intedges \\ \sigma=K|L}}
	\left[ (-F^{(n)}_{K,\sigma})^{+}c^{(n)}_K
	- (-F^{(n)}_{K,\sigma})^{-}c^{(n)}_L\right]\varphi^{(n)}_K\\
	+\int_0^T\int_\O\left[ q^P(x,t)\Pi_\mesh c(x,t)\Pi_\mesh\varphi(x,t) 
	- q^I(x,t)\hat{c}(x,t)\Pi_\mesh\varphi(x,t)\right]\ud x\ud t = 0,
	\end{gathered}\right\} \label{eq:scheme}
\end{equation}
where $(-F^{(n)}_{K,\sigma})^{+}$ and $(-F^{(n)}_{K,\sigma})^{-}$ denote the positive and
negative parts of $-F^{(n)}_{K,\sigma}$, respectively.

The following convergence is known:
\begin{theorem}[Chainais-Hillairet--Droniou \cite{cd07}] \label{th:cd}
Take a sequence $(\mesh_m)_{m\in\NN}$ of admissible meshes of $\O\times(0,T)$ satisfying
the appropriate regularity hypotheses and such that
$\size{\mesh_m}\to0$ and $\delta_{\mesh_{m}}\to 0$ as $m\to\infty$.
Then, up to a subsequence, 
\begin{enumerate}[(i)]
\item $\Pi_{\mesh_m}c\to\cbar$, a.e. on $\O\times(0,T)$,
weakly-$\ast$ in $\Leb{\infty}{2}$ and strongly in $\Leb{p}{q}$ for all $p<\infty$ and all $q<2$;
\item $\nabla_{\mesh_m}c\weakto\nabla\cbar$ weakly in $L^2(\O\times(0,T))^d$, and
\item $\U_m \to \Ubar$, weakly-$\ast$ in $\Leb[d]{\infty}{2}$ 
and strongly in $L^2(\O\times(0,T))^d$.
\end{enumerate}
\end{theorem}
The purpose of this note is to demonstrate that (i) holds with $p=\infty$ and $q=2$.

\section{Uniform temporal convergence} \label{sec:unifconv}

In the following estimate, we employ the dual seminorm
\begin{equation*}
\left|v\right|_{\star,\mesh} = \sup\left\{\int_\O\Pi_\mesh v(x)\Pi_\mesh w(x)\ud x : 
w\in\unknowns, \, \norm{\nabla_\mesh w}_{L^{2d}(\O)^d}=1\right\}.
\end{equation*}
\begin{lemma}[Discrete time derivative estimate] \label{lem:dtcest}
Let $\mesh$ be an admissible mesh of $\O\times(0,T)$ and take a solution
$(c^{(n)})_{n=0,\ldots,N}\subset\unknowns$, $(F^{(n)})_{n=1,\ldots,N}\subset\fluxes$ to \eqref{eq:scheme}. 
Then there exists $\ctel{time}>0$, not depending upon the mesh, such that
\begin{equation}\label{est:dtc}
\int_0^T \left| \delta_\mesh c(t)\right|^4_{\star, \mesh}\ud t \leq \cter{time}.
\end{equation}
\end{lemma}
\begin{sketch}
Let $w\in\unknowns$ and $t\in(0,T)$. For $\zeta\in\RR^d$, denote by $\Dsq{\cdot}{\zeta}$ the square root of the 
positive-definite matrix $\D{\cdot}{\zeta}$. From \eqref{hyp:D}, the estimates
\cite[Propositions 3.1, 3.2]{cd07} on $\U$ and $\nabla_\mesh c$ and 
the bound $|\Dsq{\cdot}{\U}|\leq\Lambda_\DD^{1/2}(1+ |\U|)^{1/2}$, one can show that
$\int_\O\D{\cdot}{\U}\nabla_\mesh c(t)\cdot\nabla_\mesh w
\leq\ctel{}\norm{\nabla_\mesh w}_{L^4(\O)^d}$.
Use the conservativity $F^{(n)}_{K,\sigma}+F^{(n)}_{L,\sigma}=0$ ($\sigma=K|L\in\intedges$) 
to gather by edges and write
\begin{multline*}
\sum_{K\in\cells}\sum_{\substack{\sigma\in\Kedges\cap\intedges \\ \sigma=K|L}}
\left[ (-F^{(n)}_{K,\sigma})^{+}c^{(n)}_K
- (-F^{(n)}_{K,\sigma})^{-}c^{(n)}_L\right]w^{(n)}_K\\
=\sum_{\sigma=K|L\in\intedges}(-F^{(n)}_{K,\sigma})^{+}c^{(n)}_K (w^{(n)}_K - w^{(n)}_L)
+ \sum_{\sigma=K|L\in\intedges}(-F^{(n)}_{L,\sigma})^{+}c^{(n)}_L (w^{(n)}_L - w^{(n)}_K).
\end{multline*} 
These terms are symmetric in $K$ and $L$ so it suffices to estimate only one of them.
By squaring both sides of \cite[Eq. (2.22)]{degh10} and making the appropriate scalings,
one can show that 
\begin{equation*}
\sum_{\sigma=K|L\in\intedges}\frac{d_{K,\sigma}}{|\sigma|}\left|-F^{(n)}_{K,\sigma}\right|^2\leq\ctel{}.
\end{equation*}
Then Cauchy-Schwarz and H\"older give, for some $p>1$ to be determined,
\begin{align*}
&\left|\sum_{\sigma=K|L\in\intedges}(-F^{(n)}_{K,\sigma})^{+}c^{(n)}_K (w^{(n)}_K - w^{(n)}_L)\right|\\
&\leq\ctel{}\left(\sum_{K\in\cells}\sum_{\substack{\sigma\in\Kedges\cap\intedges\\\sigma=K|L}}
|\sigma|d_{K,\sigma}\left|\frac{w^{(n)}_K-w^{(n)}_L}{d_{K,\sigma}}\right|^{2p'}\right)^{\frac{1}{2p'}}
\left(\sum_{K\in\cells}\sum_{\sigma\in\Kedges}
(|\sigma|d_{K,\sigma})|c^{(n)}_K|^{2p}\right)^{\frac{1}{2p}}\\
&\leq\ctel{}\norm{\nabla_\mesh w}_{L^{2p'}(\O)^d}\norm{\Pi_\mesh c(t)}_{L^{2p}(\O)},
\end{align*}
the last inequality coming from \cite[Lemma 4.2]{egh10} and 
the identity $\sum_{\sigma\in\Kedges}|\sigma|d_{K,\sigma}=d|K|$. 
Applying a discrete Sobolev inequality \cite[Lemma 5.3]{egh10} to 
the estimates \cite[Proposition 3.2]{cd07}, we can interpolate between the spaces
$\Leb{\infty}{2}$ and $\Leb{2}{2^{\star}}$ to ensure that $\Pi_\mesh c$ is bounded
in the $\Leb{4}{\frac{2d}{d-1}}$ norm. We therefore set $p=\frac{d}{d-1}>1$, which
gives $p'=d$, thereby justifying the choice of Lebesgue exponent in the dual seminorm 
definition above.

Denote by $k$ the index such that $t\in[t^{(k-1)},t^{(k)})$. In \eqref{eq:scheme},
take $\varphi=(\varphi^{(n)})_{n=1,\ldots,N}\subset\unknowns$ satisfying $\varphi^{(k)}=w$ and 
$\varphi^{(n)}=0$ for $n\neq k$ to see that 
\begin{multline*}
\int_\O\Pi_\mesh\delta_\mesh c(t)\Pi_\mesh w
\leq\ctel{}\norm{\nabla_\mesh w}_{L^{2d}(\O)^d}\bigg(
1 + \norm{\Pi_\mesh c(t)}_{L^{\frac{2d}{d-1}}(\O)}\\ 
+ \norm{q^P}_{\Leb{\infty}{r}}\norm{\Pi_\mesh c(t)}_{L^2(\O)}
+ \norm{q^I}_{\Leb{\infty}{2}}\bigg). 
\end{multline*}
Then
\begin{equation*}
\int_0^T\left| \delta_\mesh c(t)\right|^4_{\star,\mesh}\ud t
\leq\ctel{}\norm{\Pi_\mesh c(t)}_{\Leb{4}{\frac{2d}{d-1}}}\leq\cter{time}.
\qedhere
\end{equation*}
\end{sketch}
The key ideas for the following result are due to \cite{dey15}.
\begin{theorem}[Uniform temporal convergence of concentration]
Assume the same hypotheses as Theorem \ref{th:cd}. If, for all 
$T_0\in[0,T]$, $\cbar$ and $\Ubar$ satisfy
the energy identity
\begin{multline} \label{eq:contenergyeq}
\frac{1}{2}\int_\O \Phi\cbar(T_0)^2
= \frac{1}{2}\int_\O \Phi c_0^2 
+ \int_0^{T_0}\int_\O \cbar\hat{c}q^I 
- \frac{1}{2}\int_0^{T_0}\int_\O \cbar^2(q^I + q^P)
- \int_0^{T_0}\int_\O\D{\cdot}{\Ubar}\nabla\cbar\cdot\nabla\cbar,
\end{multline}
then, up to a subsequence, $\Pi_{\mesh_m}c\to\cbar$ in $\Leb{\infty}{2}$.
\end{theorem}
\begin{remark}
The identity \eqref{eq:contenergyeq} seems natural, but certainly not obvious. 
If $\DD$ is uniformly bounded, its proof is straightforward; see for example
the calculations in \cite[Proposition 3.1]{dt14}. 
\end{remark}
\begin{sketch}
Fix $T_0\in[0,T]$ and take a sequence $(T_m)_{m\in\NN}\subset[0,T]$ with $T_m\to T_0$
as $m\to\infty$. Denote by $k_m\in[1,N]$ the index such that $T_0\in[t^{(k_m-1)},t^{(k_m)})$.
Apply the uniform-in-time, weak-in-space discrete Aubin--Simon 
theorem \cite[Theorem 3.1]{dey15} with estimates \cite[Proposition 3.2]{cd07} 
and \eqref{est:dtc} to obtain 
$\Pi_{\mesh_m} c\to\cbar$ in $L^\infty(0,T;L^{2}(\O)\weak)$. This gives 
$\sqrt{\Phi}\Pi_{\mesh_m}c(T_m)\weakto\sqrt{\Phi}\cbar(T_0)$ weakly in $L^2(\O)$ and hence
\begin{equation} \label{eq:liminf}
\liminf_{m\to\infty}\int_\O\Phi\left(\Pi_{\mesh_m}c(T_m)\right)^2\geq\int_\O\Phi\left(\cbar(T_0)\right)^2.
\end{equation}
Take $\varphi=(c^{(1)},\ldots,c^{(k_m)},0,\ldots,0)\subset\unknowns$ in \eqref{eq:scheme}
and follow the calculations in \cite[Proposition 3.2]{cd07}
to obtain
\begin{multline*}
\frac{1}{2}\int_\O\Phi\left(\Pi_{\mesh_m} c(T_m)\right)^2 
+ \int_0^{T_m}\int_\O\D{\cdot}{\U_m}\nabla_{\mesh_m}c\cdot\nabla_{\mesh_m}c\\
+ \frac{1}{2}\int_0^{T_m}\int_\O (\Pi_{\mesh_m}c)^2(q^I + q^P)
\leq \frac{1}{2}\int_\O\Phi (\cI_{\mesh_m}c_0)^2
+ \int_0^{t^{(k_m)}}\int_\O \Pi_{\mesh_m}c\hat{c}q^I. 
\end{multline*}
Using the fact that $\cI_{\mesh_m}c_0\to c_0$ in $L^2(\O)$, take the limit superior as $m\to\infty$:
\begin{multline} \label{eq:terms}
\frac{1}{2}\limsup_{m\to\infty}\int_\O\Phi(\Pi_{\mesh_m}c(T_m))^2
\leq \frac{1}{2}\int_\O\Phi c^2_0
+ \limsup_{m\to\infty}\int_0^{t^{(k_m)}}\int_\O\Pi_{\mesh_m}c\hat{c}q^I\\
-\frac{1}{2}\liminf_{m\to\infty}\int_0^{T_m}\int_\O(\Pi_{\mesh_m}c)^2(q^I + q^P)
-\liminf_{m\to\infty}\int_0^{T_m}\int_\O\D{\cdot}{\U_m}\nabla_{\mesh_m}c\cdot\nabla_{\mesh_m}c\\
=: \frac{1}{2}\int_\O\Phi c^2_0 + \limsup_{m\to\infty}\cS_1^{(m)}
-\frac{1}{2}\liminf_{m\to\infty}\cS_2^{(m)} - \liminf_{m\to\infty}\cS_3^{(m)}.
\end{multline}
Note that the $t^{(k_m)}$ such that $T_m\in[t^{(k_m-1)},t^{(k_m)})$
converges to $T_0$ as $m\to\infty$. From Theorem \ref{th:cd} and Fatou's lemma,
\begin{equation*}
\limsup_{m\to\infty}\cS_1^{(m)} = \int_0^{T_0}\int_\O\cbar\hat{c}q^I \quad\mbox{and}\quad
\liminf_{m\to\infty}\cS_2^{(m)}\geq\int_0^{T_0}\int_\O\cbar^2(q^I + q^P).
\end{equation*}
A similar argument to \cite[Remark 3.2]{dt14} shows that
$\Dsq{\cdot}{\U_m}\nabla_{\mesh_m}c\weakto\Dsq{\cdot}{\Ubar}\nabla\cbar$ weakly in
$\Leb[d]{2}{2}$. Dominated convergence gives 
$\ch_{[0,T_m]}\Dsq{\cdot}{\Ubar}\nabla\cbar\to\ch_{[0,T_0]}\Dsq{\cdot}{\Ubar}\nabla\cbar$
in $\Leb[d]{2}{2}$. Thus
\begin{multline*}
\int_0^{T_0}\int_\O\D{\cdot}{\Ubar}\nabla\cbar\cdot\nabla\cbar
= \lim_{m\to\infty}\int_0^{T_m}\int_\O\left(\Dsq{\cdot}{\Ubar}\nabla_\mesh\cbar\right)
\left(\Dsq{\cdot}{\U_m}\nabla_{\mesh_m}c\right) \\
\leq\norm{\Dsq{\cdot}{\Ubar}\nabla\cbar}_{\Leb[d]{2}{2}}
\liminf_{m\to\infty}\norm{\ch_{[0,T_m]}\Dsq{\cdot}{\U_m}\nabla_{\mesh_m}c}_{\Leb[d]{2}{2}},
\end{multline*}
and so 
\begin{equation*}
\liminf_{m\to\infty}\cS_3^{(m)}\geq\int_0^{T_0}\int_\O\D{\cdot}{\Ubar}\nabla\cbar\cdot\nabla\cbar.
\end{equation*}
Collecting these convergences, we see that the right-hand sides of \eqref{eq:terms} and 
\eqref{eq:contenergyeq} agree, giving
\begin{equation} \label{eq:limsup}
\limsup_{m\to\infty}\int_\O\Phi(\Pi_{\mesh_m}c(T_m))^2\leq\int_\O\Phi\cbar(T_0)^2.
\end{equation}
Comparing \eqref{eq:liminf} and \eqref{eq:limsup} shows that 
$\lim_{m\to\infty}\norm{\sqrt{\Phi}\Pi_{\mesh_m}c(T_m)}^2_{L^2(\O)}=\norm{\sqrt{\Phi}\cbar(T_0)}^2_{L^2(\O)}$,
which, thanks to the weak-$L^2(\O)$ convergence established earlier, gives
$\sqrt{\Phi}\Pi_{\mesh_m}c(T_m)\to\sqrt{\Phi}\cbar(T_0)$ strongly in $L^2(\O)$.
Apply the characterisation \cite[Lemma 6.4]{dey15} of uniform convergence and the
uniform positivity of $\Phi$ to conclude the proof.
\end{sketch}


\bibliography{peaceman-hmm-submitted}{}
\bibliographystyle{plain}

\end{document}